\documentclass[a4paper]{article}
\usepackage{amsmath}
\usepackage{amsfonts}
\usepackage{amssymb}
\usepackage{amsthm}
\usepackage{graphicx}
\usepackage{color}
\usepackage{enumitem}
\usepackage{bm}
\usepackage{mathtools}
\usepackage{enumitem}
\usepackage{leftidx}
\usepackage{isomath}
\usepackage[colorlinks=true, urlcolor=cyan, citecolor=magenta, breaklinks=true]{hyperref}
\usepackage[hyphenbreaks]{breakurl}
\newcommand{\remove}[1] {}

\newtheorem*{summary*}{Summary of results}

\newcommand{\eat}[1]{}

\author {%Author name(s) and affiliation removed\remove{
{ Lefteris Kirousis and Georgios Kontogeorgiou}\\[0.1cm] Department of Mathematics\\National and Kapodistrian University of Athens\\
 \href{mailto:lkirousis@math.uoa.gr}{lkirousis@math.uoa.gr},\ \href{mailto:georgioskontogeorgiou97@gmail.com}  {georgioskontogeorgiou97@gmail.com}}
\title{The {\em  probl\`{e}me des   m\'{e}nages}   revisited}

\begin{document}

\maketitle
\begin{abstract}
We present an alternative proof to the Touchard-Kaplansky formula for the {\em  probl\`{e}me des   m\'{e}nages}, which, we believe, is simpler than the extant ones and is in the spirit of the elegant original proof  by Kaplansky (1943). About the  latter proof,  Bogart and Doyle (1986) argued  that despite its cleverness,  suffered from opting to give  precedence to one of the genders for the couples involved (Bogart and Doyle supplied an elegant proof that avoided such gender-dependent bias).   
\end{abstract}

%\section{Introduction} \label{intro:sec}
\bigskip\noindent
The {\em  probl\`{e}me des   m\'{e}nages} (the couples problem) asks to count the number of ways  that $n$ man-woman couples can sit around a circular table so that no one  sits next to her partner or someone of the same gender. 

The problem was first stated in an equivalent but different form by Tait \cite[p. 159]{tait1877knots} in the framework of knot theory: 
\begin{quote}
``How many arrangements are there of n letters, when A cannot be in the first or second place, B not in the second or third, \&c."	
\end{quote}
It is easy to see that the m\'{e}nage problem is equivalent  to Tait's arrangement problem:  first sit around the table, in $2n!$ ways,  the men or alternatively the women, leaving an empty space between any two consecutive of them,  then count the number of ways to sit the members of the other gender.

Recurrence relations that compute the answer to Tait's question were  soon given by Cayley and Muir \cite{Cayley1878problem, muir1878professor, Cayley1878mr, muir1882additional}. Almost fifteen years later, Lucas \cite{lucas1891théorie}, evidently  unaware of the work of Tait, Cayley and Muir, posed the problem in the formulation of husbands and wives, named it  {\em probl\`{e}me des   m\'{e}nages} and supplied a recurrence relation already described by Cayley and Muir. But it was not earlier than  fourty-three years later  that an explicit formula for Tait's problem was given by Touchard  \cite{touchard1934probleme}, alas without a proof (for  historic accounts, see Kaplansky and Riordan \cite{kaplansky1946probleme} and Dutka \cite{dutka1986probleme}).

\begin{figure}[h]
\centering
\includegraphics{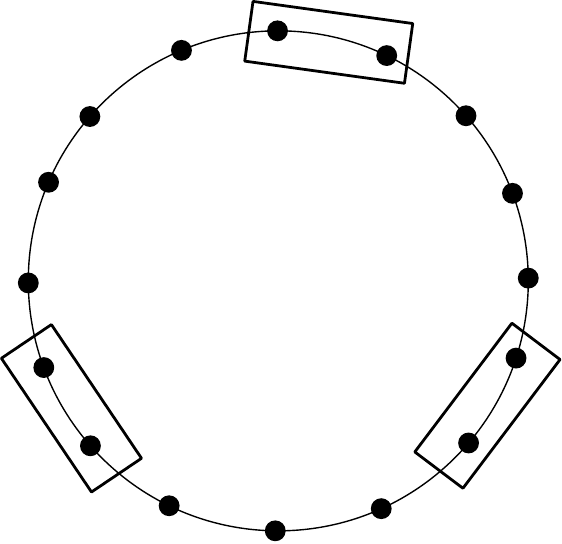}	
\caption{A cycle with 16 places and 3 non-overlapping dominos}
\label{fig:dominos}
\end{figure}

The first proof for Touchard's  explicit formula  was given nine years later  (sixty-five years later than Tait's question) by Kaplansky \cite{kaplansky1943solution}.  Specifically he showed that the number of permutations of $\{1, \ldots, n\}$ that differ in all places  with both the identity permutation and the  circular permutation $(1 \rightarrow 2 \rightarrow \cdots \rightarrow n \rightarrow 1)$ is

\begin{equation}\label{touchard}\sum_{r=0}^n (-1)^r\frac{2n}{2n-r} {{2n-r}\choose{r}} (n-r)!\end{equation}
For the proof, he used the principle of inclusion and exclusion  in a way that in \cite{bogart1986non} is characterized as ``simple but not straightforward".  As for the solution for the formulation in terms  of husbands and wives, Kaplansky and Riordan \cite{kaplansky1946probleme} in a later exposition wrote:
\begin{quote} ``We begin by fixing the positions of husbands and wives, say wives, for courtesy's sake."
 \end{quote}
Forty-three years later (more than a century later than Tait), Bogart and Doyle  in a much referenced paper in the {\em American Mathematical Monthly} \cite{bogart1986non}   gave the first proof of the explicit formula for  the m\'enage problem, {\em not} starting from a reduction to Tait's problem. They claimed that:
\begin{quote}
``Seating the ladies first ``reduces" the m\'{e}nage problem to a problem of permutations with restricted position. Unfortunately, this new problem is more difficult than the problem we began with, as we may judge from the cleverness of Kaplansky's solution."
\end{quote}
And they added as a conclusion:
\begin{quote}
``5. Conclusion. It appears that it was only the tradition of seating the ladies first that made the
m\'{e}nage problem seem in any way difficult. We may speculate that, were it not for this tradition,
it would not have taken half a century to discover Touchard's formula for $M_n$. Of all the ways in which sexism has held back the advance of mathematics, this may well be the most peculiar. (But see Execise 2.)"	
\end{quote}
See also {\em The New York Times} article \cite{thenytimes1986}. 
The  method Bogart and Doyle used  is indeed  clever and  simple. They started by counting the number of ways to place $n$ non-overlapping dominos on a cycle with $2n$ positions (see Figure \ref{fig:dominos}). This can be done in $d_r$ wqys, where:
 \begin{equation}\label{dominos} d_r = \frac{2n}{2n-r} {{2n-r}\choose{r}}.\end{equation} 
This calculation, as Bogart and Doyle write,  is a routine combinatorial problem, which  they leave as an exercise and also give a reference (we give below an easy proof for completeness). Then  they used  the principle of inclusion and exclusion by first  counting  in how many ways the  members of the couples can  be seated, so that  no two members of the same gender are adjacent,  and  at least $r$ of the couples  occupy each a single domino.  However, this counting  takes some effort. Here is how to avoid even this, perhaps at the expense of political correctness!
\paragraph{Alternative proof using both dominos and the ``sexist" reduction to Tait's problem} To avoid gender-based language, let us say that we are given $n$ couples of the same, for each couple,  letter,   in $X$-font and in $Y$-font.  We  start by fixing the positions of the  $X$-letters.  
There are $2n!$ ways to do this. Then, aiming at using the principle of inclusion and exclusion again,  we place, in $d_r$ ways, $r$ {\em transparent} non-overlapping dominos on the cycle. Now comes what apparently escaped Bogart and Doyle, as  they did not start with all $X$-letters positioned first:  there are $(n-r)!$ ways to allocate the $n$ $Y$-letters so that at least $r$ of them are seated next to their counterpart; indeed $r$ {\em uniquely determined} $Y$-letters will take the free spaces in respective $r$  dominos, whose  $X$-letter is in common view (dominos are transparent); the rest will be arbitrarily assigned to the remaining empty spaces of the cycle. So all in all, by the principle of inclusion and exclusion, the answer to the {\em  probl\`{e}me des   m\'{e}nages} is 
\begin{equation*}2n!\sum_{r=0}^n (-1)^r\frac{2n}{2n-r} {{2n-r}\choose{r}} (n-r)! \tag*{$\Box$}\end{equation*}
\paragraph{Proof of equation \eqref{dominos}}  Pick a starting place on the cycle: $2n$ ways; place the identical dominos with undetermined number of empty spaces in each of  the $r$  
arcs between any two consecutive of them, these arcs should be filled with $2n-2r$ empty spaces; determine the number of empty spaces in each arc between two consecutive dominos by throwing $2n-2r$ identical balls into $r$ distinguished bins: ${2n-r+1}\choose{r-1}$ ways; divide by $r$ since in the counting so far, each of the  circularly placed $r$ identical dominos counted separately as the first one, to finally get: 
\begin{equation*}d_r = \frac{2n}{r} {{2n-r-1}\choose{r-1}} = \frac{2n}{2n-r} {{2n-r}\choose{r}}. \tag*{$\Box$}\end{equation*} 

\paragraph{Conclusion} It is quite unfortunate that putting  some non-mathematical remarks of a 1943 proof for  a 19th century problem  under a contemporary social lens  hindered  the very clever  idea of using dominos to show its full simplifying power by assuming that underneath the dominos lie the prioritized   $X$-letters.

%\bibliographystyle{plain}
%\bibliography{menage}

\end{document}